\documentclass{amsart}

\newcommand{\D}{d}

\newtheorem{theorem}{Theorem}[section]
\newtheorem{proposition}[theorem]{Proposition}
\newtheorem{lemma}[theorem]{Lemma}

\theoremstyle{definition}
\newtheorem{definition}[theorem]{Definition}

\theoremstyle{remark}

\begin{document}

\title{The Geography of Non-formal Manifolds}
\author{Marisa Fern\'andez\and Vicente Mu\~noz}

\maketitle

\vskip1cm

\begin{abstract}
We show that there exist non-formal compact oriented manifolds of
dimension $n$ and with first Betti number $b_1=b\geq 0$ if and
only if $n\geq 3$ and $b\geq 2$, or $n\geq (7-2b)$ and $0\leq
b\leq 2$. Moreover, we present explicit examples for each one of
these cases.
\end{abstract}

\bigskip

\section{Introduction} \label{intro}

Simply connected compact manifolds of dimension less than or equal
to $6$ are formal \cite{NM,Mi,FM}. A method to construct
non-formal simply connected compact manifolds of any dimension
$n\geq 7$ was given by the authors in \cite{FM2}. An alternative
method is given in \cite{Ru} (see also \cite{O} for an example in
dimension $7$). A natural question to ask is whether there are
examples of non-formal compact manifolds of any dimension whose
first Betti number $b_1=b \geq 0$ is arbitrary.

We consider the following problem on the \emph{geography} of
manifolds:

\noindent For which pairs $(n,b)$ with $n\geq 1$ and $b\geq 0$ are
there compact oriented manifolds of dimension $n$ and with $b_1=b$
which are non-formal?

Note that we can restrict to just considering connected manifolds.
In this paper we solve completely this problem by proving the
following main result.

\begin{theorem} \label{thm:main}
 There are compact oriented $n$-dimensional manifolds with $b_1=b$
 which are non-formal if and only if $n\geq 3$
 and $b\geq 2$, or $n\geq (7-2b)$ and $0\leq b\leq 2$.
\end{theorem}

In the case of a simply connected manifold $M$, formality for $M$
is equivalent to saying that its real homotopy type is determined
by its real cohomology algebra. In the non-simply connected case,
things are a little bit more complicated. If $M$ is nilpotent,
i.e., $\pi_1(M)$ is nilpotent and it acts nilpotently on
$\pi_i(M)$ for $i\geq2$, then formality means again that the real
homotopy type is determined by the real cohomology algebra. In
general, we shall say that $M$ is formal if the minimal model of
the manifold (which is, by definition, the minimal model of the
algebra of differential forms $\Omega^*(M)$) is determined by the
real cohomology algebra (see Sect.\ \ref{definitions} for precise
definitions). Note that there are alternative (and non-equivalent)
definitions of formality in the non-nilpotent situation (see
\cite{Tanre}). This punctualization is important because the
non-formal manifolds that we construct in Sect.\ \ref{nonformal}
are necessarily not nilpotent (see Sect.\ \ref{remarks}).

In the following table, the big dots mark the pairs $(n,b_1)$ for
which all manifolds of dimension $n$ and first Betti number $b_1$
are formal. For any of the small dots, there are examples of
non-formal manifolds.

\begin{table}
\centering \caption{Geography of non-formal manifolds}
\label{tab:1}
 $
 \begin{array}{r|cccc}
 n\geq7 & \cdot & \cdot &\cdot & \cdots \\
 n=6 & \bullet & \cdot & \cdot & \cdots \\
 n=5 & \bullet &\cdot & \cdot & \cdots\\
 n=4 & \bullet & \bullet & \cdot & \cdots \\
 n=3 & \bullet & \bullet & \cdot & \cdots \\
 n=2 & \bullet & \bullet & \bullet & \bullet \\
 \hline  & \quad b_1=0 \quad& b_1=1&  \quad b_1=2 \quad &  b_1\geq 3\quad
 \end{array}
 $
\end{table}

To prove Thm.\ \ref{thm:main} we need to do two things. On the one
hand, we need to verify that manifolds of dimension $n\leq 6$ with
$b_1=0$ and manifolds of dimension $n\leq 4$ with $b_1=1$ are
\emph{always} formal. For this we use the results of \cite{FM}. On
the other hand, we need to present examples of \emph{non-formal}
manifolds of dimension $n\geq 7$ with $b_1=0$, of dimension $n\geq
5$ with $b_1=1$ and of dimension $n\geq 3$ for any other $b_1\geq
2$. For this we use a similar method to that of \cite{FM2}. Note
that both questions for the case $b_1=0$ are already solved, so
here we have to focus on the case $b_1=1$.


\section{Minimal Models and Formality} \label{definitions}

We recall some definitions and results about minimal models
\cite{H,DGMS,Ta}. Let $(A,\D )$ be a {\it differential algebra},
that is, $A$ is a graded commutative algebra over the real
numbers, with a differential $\D $ which is a derivation, i.e. $\D
(a\cdot b) = (\D  a)\cdot b +(-1)^{\deg (a)} a\cdot (\D  b)$,
where $\deg(a)$ is the degree of $a$. Morphisms between
differential algebras are required to be degree preserving algebra
maps which commute with the differentials.

A differential algebra $(A,\D )$ is said to be {\it minimal\/} if:
\begin{enumerate}
 \item $A$ is free as an algebra, that is, $A$ is the free
 algebra $\bigwedge V$ over a graded vector space $V=\oplus V^i$, and
 \item there exists a collection of generators $\{ a_\tau,
 \tau\in I\}$, for some well ordered index set $I$, such that
 $\deg(a_\mu)\leq \deg(a_\tau)$ if $\mu < \tau$ and each $\D
 a_\tau$ is expressed in terms of preceding $a_\mu$ ($\mu<\tau$).
 This implies that $\D  a_\tau$ does not have a linear part, i.e., it
 lives in ${\bigwedge V}^{>0} \cdot {\bigwedge V}^{>0} \subset {\bigwedge V}$.
\end{enumerate}

We shall say that a minimal differential algebra $(\bigwedge V,\D
)$ is a {\it minimal model} for a connected differentiable
manifold $M$ if there exists a morphism of differential graded
algebras $\rho\colon {(\bigwedge V,\D )}\longrightarrow {(\Omega
M,\D )}$, where $\Omega M$ is the de Rham complex of differential
forms on $M$, inducing an isomorphism
 $$
 \rho^*\colon H^*(\bigwedge V)\longrightarrow H^*(\Omega M,\D )= H^*(M)
 $$
on cohomology.

If $M$ is a simply connected manifold (or, more generally, a
nilpotent space), the dual of the real homotopy vector space
$\pi_i(M)\otimes {{\mathbf{R}}}$ is isomorphic to $V^i$ for any
$i$. Halperin in~\cite{H} proved that any connected manifold $M$
has a minimal model unique up to isomorphism, regardless of its
fundamental group.

A minimal model $(\bigwedge V,\D )$ of a manifold $M$ is said to
be {\it formal}, and $M$ is said to be {\it formal}, if there is a
morphism of differential algebras $\psi\colon {(\bigwedge V,\D
)}\longrightarrow (H^*(M),\D =0)$ that induces the identity on
cohomology. An alternative way to look at this is the following:
the above property means that $(\bigwedge V,\D )$ is a minimal
model of the differential algebra $(H^*(M),0)$. Therefore $(\Omega
M,\D )$ and $(H^*(M),0)$ share their minimal model, i.e., one can
obtain the minimal model of $M$ out of its real cohomology
algebra. When $M$ is nilpotent, the minimal model encodes its real
homotopy type.

\bigskip

In order to detect non-formality, we have Massey products. Let us
recall its definition. Let $M$ be a (not necessarily simply
connected) manifold and let $a_i \in H^{p_i}(M)$, $1 \leq i\leq
3$, be three cohomology classes such that $a_1\cup a_2=0$ and
$a_2\cup a_3=0$. Take forms $\alpha_i$ in $M$ with
$a_i=[\alpha_i]$ and write $\alpha_1\wedge \alpha_2= \D \xi$,
$\alpha_2\wedge \alpha_3=\D \eta$. The Massey product of the
classes $a_i$ is defined as
  $$
  \langle a_1,a_2,a_3 \rangle  = [ \alpha_1 \wedge \eta+(-1)^{p_1+1} \xi
  \wedge \alpha_3] \in \frac{H^{p_1+p_2+p_3-1}(M)}{a_1
  \cup H^{p_2+p_3-1}(M) + H^{p_1+p_2-1}(M)\cup a_3} \; .
  $$

We have the following result, for whose proof we refer
to~\cite{DGMS,Ta,TO}.

\begin{theorem} \label{thm:criterio1}
 If $M$ has a non-trivial Massey product then $M$ is non-formal.
\end{theorem}

Therefore the existence of a non-zero Massey product is an
obstruction to the formality.

\bigskip

In order to prove formality, we extract the following notion from
\cite{FM}.

\begin{definition}\label{def:s-formal}
 Let $(\bigwedge V,\D )$ be a minimal model of a differentiable manifold
$M$. We say that $(\bigwedge V,\D )$ is \emph{$s$-formal}, or $M$
is a \emph{$s$-formal manifold} ($s\geq 0$) if for each $i\leq s$
one can get a space of generators $V^i$ of elements of degree $i$
that decomposes as a direct sum $V^i=C^i\oplus N^i$, where the
spaces $C^i$ and $N^i$ satisfy the three following conditions:
 \begin{enumerate}
 \item $\D (C^i) = 0$,
 \item the differential map $\D \colon N^i\longrightarrow \bigwedge V$ is
   injective,
 \item any closed element in the ideal $I_s=
   I(\bigoplus\limits_{i\leq s} N^i)$, generated by
   $\bigoplus\limits_{i\leq s} N^i$ in $\bigwedge
   (\bigoplus\limits_{i\leq s} V^i)$, is exact in $\bigwedge V$.
 \end{enumerate}
\end{definition}

The condition of $s$-formality is weaker than that of formality.
However we have the following positive result proved in \cite{FM}.

\begin{theorem}\label{thm:criterio2}
Let $M$ be a connected and orientable compact differentiable
manifold of dimension $2n$ or $(2n-1)$. Then $M$ is formal if and
only if is $(n-1)$-formal.
\end{theorem}

This result is very useful because it allows us to check that a
manifold $M$ is formal by looking at its $s$-stage minimal model,
that is, $\bigwedge (\bigoplus\limits_{i\leq s} V^i)$. In general,
when computing the minimal model of $M$, after we pass the middle
dimension, the number of generators starts to grow quite
dramatically. This is due to the fact that Poincar\'e duality
imposes that the Betti numbers do not grow and therefore there are
a large number of cup products in cohomology vanishing, which must
be killed in the minimal model by introducing elements in $N^i$,
for $i$ above the middle dimension. This makes Thm.\
\ref{thm:criterio2} a very useful tool for checking formality in
practice.

\section{Non-Formal Manifolds with $b_1=1$ and Dimensions $5$ and $6$}
\label{nonformal}

\subsection*{The $5$-Dimensional Example}

Let $H$ be the Heisenberg group, that is, the connected nilpotent
Lie group of dimension $3$ consisting of matrices of the form
 $$
 a=\left( \begin{array}{ccc} 1&x&z\\ 0&1&y\\ 0&0&1 \end{array}\right),
 $$
where $x,y,z \in {{\mathbf{R}}}$. Then a global system of
coordinates ${x,y,z}$ for $H$ is given by $x(a)=x$, $y(a)=y$,
$z(a)=z$, and a standard calculation shows that a basis for the
left invariant $1$-forms on $H$ consists of $\{\D  x, \D  y, \D
z-x\,\D  y\}$. Let $\Gamma$ be the discrete subgroup of $H$
consisting of matrices whose entries are integer numbers. So the
quotient space $N =\Gamma{\backslash}H$ is a compact
$3$-dimensional nilmanifold. Hence the forms $\D  x$, $\D  y$, $\D
z-x\,\D  y$ descend to $1$-forms $\alpha$, $\beta$, $\gamma$ on
$N$ and
 $$
 \D  \alpha=\D  \beta=0, \quad  \D  \gamma=-\alpha \wedge \beta\; .
 $$
The non-formality of $N$ is detected by a non-zero triple Massey
product
 $$
 \langle [\beta],[\alpha],[\alpha]\rangle  = [\alpha \wedge \gamma] \in
 \frac{H^2(N)}{[\beta]\cup H^1(N) + H^1(N) \cup [\alpha]}=H^2(N)\; .
 $$

Now let us consider the $5$-dimensional manifold $X=N \times
{\mathbf{T}}^2$, where
${\mathbf{T}}^2={\mathbf{R}}^2/{\mathbf{Z}}^2$. The coordinates of
${\mathbf{R}}^2$ will be denoted $x_1,x_2$. So $\{\D  x_1,\D
x_2\}$ defines a basis $\{\delta_1,\delta_2\}$ for the $1$-forms
on ${\mathbf{T}}^2$. We get a non-zero triple Massey product as
follows
 \begin{equation} \label{eqn:massey}
 \langle [\beta\wedge \delta_1],[\alpha],[\alpha]\rangle  =[\gamma \wedge
\alpha\wedge \delta_1].
 \end{equation}

Our aim now is to kill the fundamental group of $X$ by performing
a suitable surgery construction, in order to obtain a manifold
with $b_1=1$.

The projection $p(x,y,z)=(x,y)$ describes $N$ as a fiber bundle
$p:N\to {\mathbf{T}}^2$ with fiber ${\mathbf{S}}^1$. Actually, $N$
is the total space of the unit circle bundle of the line bundle of
degree $1$ over the $2$-torus. The fundamental group of $N$ is
therefore
  \begin{equation}\label{eqn:pi}
  \pi_1(N) \cong \Gamma = \langle \lambda_1,\lambda_2,\lambda_3 \,|\,
[\lambda_1,\lambda_2]=\lambda_3,
    \hbox{$\lambda_3$ central} \rangle  \; ,
  \end{equation}
where $\lambda_3$ corresponds to the fiber. The fundamental group
of $X=N\times{\mathbf{T}}^2$ is
 \begin{equation}\label{eqn:pi1}
 \pi_1(X)=\pi_1(N) \oplus {\mathbf{Z}}^2\; .
 \end{equation}

Consider the following submanifolds embedded in $X$:
\begin{eqnarray*}
 T_1&=& p^{-1}(\{0\}\times  {\mathbf{S}}^1 ) \times \{0\} \times \{0\}\; , \\
 T_2&=& \{\xi \} \times {\mathbf{S}}^1 \times {\mathbf{S}}^1\; ,
\end{eqnarray*}
with $\xi$ a point in $N$. These are $2$-dimensional tori with
trivial normal bundle.

Consider now another $5$-manifold $Y$ with an embedded
$2$-dimensional torus $T$ with trivial normal bundle. Then we may
perform the \emph{fiber connected sum} of $X$ and $Y$ identifying
$T_1$ and $T$, denoted $X \#_{T_1=T} Y$, in the following way:
take (open) tubular neighborhoods $\nu_1 \subset X$ and $\nu
\subset Y$ of $T_1$ and $T$ respectively; then $\partial \nu_1
\cong {{\mathbf{T}}^2}\times {\mathbf{S}}^2$ and $\partial \nu
\cong {{\mathbf{T}}^2}\times {\mathbf{S}}^2$; take an orientation
reversing diffeomorphism $\phi:\partial \nu_1
{\stackrel{\simeq}{\rightarrow}}
\partial \nu$; the fiber connected sum is defined to be the
(oriented) manifold obtained by gluing $X-\nu_1$ and $Y-\nu$ along
their boundaries by the diffeomorphism $\phi$. In general, the
resulting manifold depends on the identification $\phi$, but this
will not be relevant for our purposes.

\begin{lemma}\label{lem:pi1}
  Suppose $Y$ is simply connected. Then the fundamental group of
  $X \#_{T_1=T} Y$
  is the quotient of $\pi_1(X)$ by the image of $\pi_1(T_1)$.
\end{lemma}

\begin{proof}
 Since the codimension of $T_1$ is bigger than or equal to $3$, we
 have that
 $\pi_1(X-\nu_1)=\pi_1(X-T_1)$ is isomorphic to $\pi_1(X)$. The
 Seifert--Van Kampen theorem establishes that $\pi_1(X \#_{T_1=T} Y)$ is
 the amalgamated sum of $\pi_1(X-\nu_1)=\pi_1(X)$ and
 $\pi_1(Y-\nu)=\pi_1(Y)=1$ over the image of $\pi_1(\partial
 \nu_1)=\pi_1(T_1\times{\mathbf{S}}^2)=\pi_1(T_1)$, as required.
\end{proof}

We shall take for $Y$ the $5$-sphere ${\mathbf{S}}^5$. We embed a
$2$-dimensional torus ${\mathbf{T}}^2$ in ${\mathbf{R}}^5$. This
torus has a trivial normal bundle since its tangent bundle is
trivial (being parallelizable) and the tangent bundle of
${\mathbf{R}}^5$ is also trivial. After compactifying
${\mathbf{R}}^5$ by one point we get a $2$-dimensional torus
$T\subset {\mathbf{S}}^5$ with trivial normal bundle.

In the same way, we may consider another copy of the
$2$-dimensional torus $T\subset {\mathbf{S}}^5$ and perform the
fiber connected sum of $X$ and ${\mathbf{S}}^5$ identifying $T_2$
and $T$. We may do both fiber connected sums along $T_1$ and $T_2$
simultaneously, since $T_1$ and $T_2$ are disjoint. Call
 \begin{equation}\label{eqn:pi1'}
  M= X \#_{T_1=T} {\mathbf{S}}^5 \#_{T_2=T} {\mathbf{S}}^5
 \end{equation}
the resulting manifold. By Lem.\ \ref{lem:pi1}, $\pi_1(M)$ is the
quotient of $\pi_1(X)$ by the images of $\pi_1(T_1)$ and
$\pi_1(T_2)$. This kills the ${\mathbf{Z}}^2$ summand in
(\ref{eqn:pi1}) and it also kills $\lambda_2$ and $\lambda_3$ in
(\ref{eqn:pi}). Therefore $\pi_1(M)=\langle \lambda_1\rangle \cong
{\mathbf{Z}}$, i.e.,\, $b_1(M)=1$.

\bigskip

Our goal is now to prove that $M$ is non-formal. We shall do this
by proving the non-vanishing of a suitable triple Massey product.
More specifically, let us prove that the Massey product
(\ref{eqn:massey}) survives to $M$.

For this, let us describe geometrically the cohomology classes
$[\alpha\wedge \delta_1]$ and $[\beta]$. Consider the following
submanifolds of $X$:
\begin{eqnarray*}
 B_1&=&p^{-1}({\mathbf{S}}^1\times\{a_2\})\times\{b_1\}\times {\mathbf{S}}^1\;
,\\
 B_2&=&p^{-1}(\{a_1\}\times{\mathbf{S}}^1)\times {\mathbf{S}}^1
\times{\mathbf{S}}^1\; ,
\end{eqnarray*}
where the $a_i$ and $b_i$ are generic points of ${\mathbf{S}}^1$.
It is easy to check that $B_i\cap T_j=\emptyset$ for all $i$ and
$j$. So $B_i$ may be also considered as submanifolds of $M$. Let
$\eta_i$ be the $2$-forms representing the Poincar\'{e} dual to $B_i$
in $X$. By \cite{BoTu}, $\eta_i$ can be taken supported in a small
tubular neighborhood of $B_i$. Therefore the support of $\eta_i$
lies inside $X-T_1-T_2$, so we also have naturally $\eta_i \in
\Omega^2(M)$. Note that in $X$ we have clearly that
$[\eta_1]=[\beta \wedge e_1]$ and $[\eta_2]= [\alpha]$, where
$e_1$ is (the pull-back of) a differential $1$-form on
${\mathbf{S}}^1$ (considered as the first of the two circle
factors in $X=N\times{\mathbf{S}}^1\times{\mathbf{S}}^1$)
cohomologous to $\delta_1$ and supported in a neighborhood of
$b_1\in {\mathbf{S}}^1$. Thus $[\eta_1]=[\beta \wedge \delta_1]$
in $X$.

\begin{lemma} \label{Massey}
 The triple Massey product $\langle [\eta_1],[\eta_2],[\eta_2]\rangle$ is
 well-defined on $M$ and equals to $[\gamma\wedge\alpha\wedge e_1]$.
\end{lemma}

\begin{proof}
Let $\alpha'$ be the pull-back to $N$ of the $1$-form supported in
a neighborhood of $a_1$ in the first factor of
${\mathbf{S}}^1\times {\mathbf{S}}^1$ under the projection $p:N\to
{\mathbf{T}}^2$. Analogously, let $\beta'$ be the pull-back to $N$
of the $1$-form supported in a neighborhood of $a_2$ in the second
factor of ${\mathbf{S}}^1\times {\mathbf{S}}^1$. Therefore
$[\alpha']=[\alpha]$ and $[\beta']=[\beta]$. Clearly
 $$
 (\alpha' \wedge e_1) \wedge \beta'  = - \D \gamma' \wedge e_1,
 $$
where $\D \gamma'=\alpha'\wedge \beta'$. It can be supposed easily
that $\gamma'$ is zero in a neighborhood of $\xi \in N$. Therefore
the support of $\gamma'\wedge e_1$ is disjoint from $T_1$ and
$T_2$. Hence $\gamma' \wedge e_1$ is well-defined as a form in
$M$. So the triple Massey product
 $$
 \langle [\eta_1],[\eta_2],[\eta_2]\rangle =[\gamma' \wedge \alpha \wedge e_1]
 $$
is well-defined in $M$.
\end{proof}

Finally let us see that this Massey product is non-zero in
 $$
 \frac{H^3(M)}{[\beta'\wedge e_1] \cup H^1(M) + H^2(M) \cup [\alpha']}\; .
 $$
Consider $B_3=p^{-1}({\mathbf{S}}^1\times \{a_3\} ) \times
{\mathbf{S}}^1 \times \{b_2\}$, for generic points $a_3,b_2$ of
${\mathbf{S}}^1$. Then the Poincar\'{e} dual of $B_3$ is defined by a
$2$-form $\beta''\wedge e_2$ supported near $B_3$, where $\beta''$
is Poincar\'{e} dual to $p^{-1}( {\mathbf{S}}^1 \times \{a_3\})$,
$[\beta'']=[\beta]$, and $e_2$ is (the pull-back of) a
differential $1$-form on ${\mathbf{S}}^1$ (considered as the
second of the two circle factors in
$X=N\times{\mathbf{S}}^1\times{\mathbf{S}}^1$) cohomologous to
$\delta_2$ and supported in a neighborhood of $b_2\in
{\mathbf{S}}^1$. Again this $2$-form can be considered as a form
in $M$. Now for any $[\varphi] \in H^1(M)$, $[\varphi']\in H^2(M)$
we have
 $$
  ([\gamma' \wedge \alpha \wedge e_1]+ [\beta'\wedge e_1 \wedge \varphi]
+[\alpha'
  \wedge \varphi']) \cdot [\beta''\wedge e_2]=1\; ,
 $$
since the first product gives $1$, the second is zero and the
third is zero because $\alpha' \wedge \beta''$ is exact in $N$ and
hence in $M$. This result and Thm.\ \ref{thm:criterio1} prove the
following

\begin{theorem} \label{5}
 The manifold $M$, defined by~(\ref{eqn:pi1'}), is a compact oriented
 non-formal $5$-manifold with $b_1=1$.
\end{theorem}

\subsection*{The $6$-Dimensional Example}

A compact oriented non simply connected and non-formal manifold $M'$ of
dimension $6$ is obtained in an analogous fashion to the
construction of the $5$-dimensional manifold $M$. We start with
$X'=N\times {\mathbf{T}}^3$ and consider the $3$-dimensional tori
with trivial normal bundle
\begin{eqnarray*}
 T_1'&=& p^{-1}(\{0\}\times {\mathbf{S}}^1 ) \times \{0\} \times \{0\}\times
{\mathbf{S}}^1\; , \\
 T_2'&=& \{\xi \} \times {\mathbf{S}}^1 \times {\mathbf{S}}^1\times
{\mathbf{S}}^1\; .
\end{eqnarray*}
Define
 \begin{equation}\label{eqn:6d}
 M'=X' \#_{T_1'=T'} {\mathbf{S}}^6 \#_{T_2'=T'}{\mathbf{S}}^6\; ,
 \end{equation}
where $T'$ is an embedded $3$-torus in ${\mathbf{S}}^6$ with
trivial normal bundle. Then $M'$ is a non-formal $6$-manifold with
$b_1=1$, which can be proved in a similar way to Thm.\ \ref{5}.

\section{Proof of Theorem \ref{thm:main}} \label{proof}

Let us first prove the affirmative results in Thm.\
\ref{thm:main}.

\begin{proposition}\label{positive}
Let $M$ be a connected, compact and orientable manifold of
dimension $n$ and first Betti number $b_1=b$.
\begin{itemize}
  \item If $n\leq 2$ then $M$ is formal.
  \item If $n\leq 6$ and $b=0$ then $M$ is formal.
  \item If $n\leq 4$ and $b=1$ then $M$ is formal.
\end{itemize}
\end{proposition}

\begin{proof}
The first item is well-known: the circle and any oriented surface
are formal. However, it follows from Thm.\ \ref{thm:criterio2}
very easily. Since $M$ is connected, $M$ is $0$-formal. Hence $M$
is formal as $n\leq 2$.

Second item follows from \cite{FM,Mi,NM}. Let us recall briefly
the proof. Since $M$ has $b_1=0$ it follows that in the minimal
model $V^1=0$. This implies that $N^2=0$ since there are no
decomposable elements of degree $3$ and hence no element of $V^2$
can kill any element of degree $3$ in the minimal model. Thus $M$
is $2$-formal and hence formal, by Thm.\ \ref{thm:criterio2},
since $n\leq 6$.

The third item is proved similarly. Since $M$ has $b_1=1$, in the
minimal model $(\bigwedge V,\D )$ we have that $V^1=C^1$ is
generated by one element $\xi$. There cannot be any element in
$N^1$ since there are no decomposable elements of degree $2$ (the
only such element is $\xi \cdot \xi = 0$). Thus $M$ is $1$-formal
and hence formal, by Thm.\ \ref{thm:criterio2}, since $n\leq 4$.
 \end{proof}

With this result, we only need to find non-formal (connected,
compact, orientable) manifolds under the conditions $n\geq
\max\{3,7-2b_1\}$ to complete the proof of Thm.\ \ref{thm:main}.

\begin{itemize}
 \item Non-formal manifolds with $n\geq 7$ and $b_1=0$ are constructed
  by the authors in \cite{FM2}. Actually those examples are simply connected.
  An alternative method is given in \cite{Ru}. Oprea \cite{O} also constructed
  examples of dimension $7$ for other purposes.
 \item Non-formal manifolds of dimensions $n=5$ or $6$ and first Betti
  number $b_1=1$. These are the manifolds $M$ and $M'$ given by
  (\ref{eqn:pi1'}) and (\ref{eqn:6d}) in Sect.\ \ref{nonformal}.
 \item Non-formal manifolds of dimension $n\geq 7$ and $b_1=1$. Take
  the non-formal $5$-dimensional manifold $M$ of
  Sect.\ \ref{nonformal} and consider $M \times {\mathbf{S}}^{n-5}$. This is
again
  non-formal (by \cite[Lem.\ 2.11]{FM})
  and has $b_1(M \times {\mathbf{S}}^{n-5})=b_1(M)=1$.
 \item Case $n=3$ and $b_1=2$. The manifold $N$ considered as the
 beginning of Sect.\ \ref{nonformal} is non-formal.
 \item Case $n=3$ and $b_1\geq 3$. Consider $N\# (b_1-2)
 ({\mathbf{S}}^1\times {\mathbf{S}}^2)$, which is non-formal because the Massey
product
 $\langle [\beta],[\alpha],[\alpha]\rangle  = [\alpha \wedge \gamma]$ is again
 defined and non-zero (as it happened for $N$).
 \item Case $n=4$ and $b_1\geq 3$. Consider $\left( N\# (b_1-3)
 ({\mathbf{S}}^1\times {\mathbf{S}}^2)\right) \times {\mathbf{S}}^1$, which is
non-formal being a
 product of a non-formal manifold with other manifold.
 \item Case $n\geq 5$ and $b_1\geq 2$. We just consider
 $\left( N\# (b_1-2) ({\mathbf{S}}^1\times {\mathbf{S}}^2)\right) \times
{\mathbf{S}}^{n-3}$.
 \item Case $n=4$ and $b_1=2$. A non-formal example can be constructed
 by a nilmanifold which is non-formal. For example (see \cite{FGG}),
 let $E$ be the total space of the ${\mathbf{S}}^1$-bundle over $N$
 with Chern class $c_1=[\beta\wedge\gamma]\in H^2(N)$. The nilmanifold
 $E$ is defined by the equations
 $$
 \D \alpha=\D \beta=0, \quad  \D \gamma=-\alpha \wedge \beta, \quad
 \D \eta=\beta\wedge \gamma\; ,
 $$
 where $\{\alpha, \beta, \gamma, \eta\}$ is a basis for the
 differential $1$-forms on $E$. Then $[\beta] \cup [\alpha]
 =[\alpha]\cup [\alpha] =0$, so that the Massey product $\langle
 [\beta],[\alpha],[\alpha]\rangle$ is well defined, and it is non-zero
 because it is represented by the cohomology class of $\gamma
 \wedge \alpha$ which is non-zero in cohomology.
\end{itemize}

\section{Final Remarks} \label{remarks}

 Note that the examples of non-formal manifolds with $b_1=1$ that we have
 constructed have abelian fundamental group, since it is isomorphic to
${\mathbf{Z}}$.
 However, these manifolds are not nilpotent.
 Actually, if a manifold $M$ with $b_1=1$ is nilpotent then $M$
 is $2$-formal. So if furthermore, the dimension is $n\leq 6$ and $M$
 is compact oriented, then it is formal.

To prove that for a nilpotent manifold $M$ with $b_1=1$ we have
that $M$ is $2$-formal, it is enough to check that $N^2=0$. This
would follow from the fact that no decomposable element of degree
$3$ (i.e., elements in $V^1\cdot V^2$) is exact. Let $\xi$ be the
generator of $V^1$ and let $a\in V^2$ be a non-zero closed
element. Suppose that $[\xi]\cup [a]=0$ and let us reach to a
contradiction. We use the following lemma of
Lalonde-McDuff-Polterovich \cite{MLP}, which has been communicated
to us by J. Oprea.

\begin{lemma} \label{LMP}
 Suppose that $\gamma\in \pi_1(M)$, $A\in \pi_2(M)$, $h\in H^1(M;{\mathbf{Z}})$
 and $\alpha\in H^2(M;{\mathbf{Z}})$, satisfy that $h(\gamma)\neq 0$ and
$\alpha(A)\neq
 0$. Then if $\alpha\cup h=0$, the action of $\gamma$ on $A$ is non-trivial.
\end{lemma}

In our case take $h=[\xi]\in H^1(M)$ (after suitable rescaling if
necessary to make it an integral class). Let $\gamma\in\pi_1(M)$
be any element with $h(\gamma)\neq 0$. Then $h(\gamma^n)\neq 0$
for any $n>0$. Now take $\alpha=[a]$ and consider any element
$A\in\pi_2(M)$ with $\alpha (A)\neq 0$ (this exists since we are
assuming that $M$ is nilpotent and in this case
$V^2=(\pi_2(M)\otimes {\mathbf{R}})^*$). Then Lem.\ \ref{LMP}
implies that $\gamma^n$ acts on $A$ non-trivially. Hence $\gamma$
acts non-nilpotently on $\pi_2(M)$, which is a contradiction.

\medskip

We would like to end up with some questions that arise naturally
once with have answered Thm.\ \ref{thm:main}.

\begin{enumerate}

\item Are there any restrictions on the Betti numbers for the
existence of non-formal manifolds? Alternatively, solve the
following \emph{geography problem}:

\noindent For which tuples $(n,b_1,\ldots,b_s)$ with $n\geq 1$,
$s=[\frac{n}2]$ and $b_i\geq 0$ is there a compact oriented
manifold $M$ of dimension $n$, with Betti numbers $b_i(M)=b_i$,
$i=1,\ldots, s$, and which is \emph{non-formal}.

\item Another alternative question is the following:

\noindent Given a finitely presented group $\Gamma$ and an integer
$n$ with $n\geq\max\{3, 2b_1(\Gamma)-7\}$, are there always
non-formal $n$-manifolds $M$ with fundamental group $\pi_1(M)\cong
\Gamma$?

\end{enumerate}

\bigskip

\noindent {\bf Acknowledgments.} Untold thanks are due to the
Organizing Committee who worked so hard to make the conference in
Lecce a success. The first author was partially
supported by CICYT (Spain) Project BFM2001-3778-C03-02 and UPV
00127.310-E-14813/2002. The second author was supported by CICYT
Project BFM2000-0024. Also both authors are partially supported by
The European Contract Human Potential Programme, Research Training
Network HPRN-CT-2000-00101.

\vspace{0.15cm}

\noindent{\sf M. Fern\'andez:} Departamento de Matem\'aticas,
Facultad de Ciencia y Tecnolog\'{\i}a, Universidad del Pa\'{\i}s Vasco,
Apartado 644, 48080 Bilbao, Spain. {\sl E-mail:}
mtpferol@lg.ehu.es

\vspace{0.15cm}

\noindent{\sf V. Mu\~noz:} Departamento de Matem\'aticas, Facultad
de Ciencias, Universidad Aut\'onoma de Madrid, 28049 Madrid,
Spain. {\sl E-mail:} vicente.munoz@uam.es

\end{document}